# ON A NONHIERARCHICAL VERSION OF THE GENERALIZED RANDOM ENERGY MODEL[1]

By Erwin Bolthausen and Nicola Kistler

*Universität Zürich*

We introduce a natural nonhierarchical version of Derrida's generalized random energy model. We prove that, in the thermodynamical limit, the free energy is the same as that of a suitably constructed GREM.

**1. Introduction and definition of the model.** The generalized random energy model (GREM for short), introduced by Derrida [1], plays an important role in spin glass theory. Originally invented as a simple model which exhibits replica symmetry breaking at various levels, it has become clear that more interesting models, like the celebrated one of Sherrington–Kirkpatrick, exhibit GREM-like behavior in the large $N$ limit. Despite the spectacular recent progress in understanding the SK-model (see [4, 5, 7]), many issues have not been clarified at all, the most prominent one being the so-called *ultrametricity*. [A metric $d$ is called an *ultrametric* if the strengthened triangle condition holds: $d(x,z) \leq \max(d(x,y), d(y,z))$. Equivalently, two balls are either disjoint or one is contained in the other.] The GREM is of limited use to investigate this because it is hierarchically organized from the start. This favorable situation allows for a complete solution, fully confirming the so-called Parisi theory (we refer the reader to the detailed study [2] where it is also pointed out that, interestingly, the emerging ultrametricity of the Gibbs measure does not necessarily coincide with the starting hierarchical organization). Yet, from the considerations on the GREM, one gets little clue on why many systems should be ultrametric in the limit. [In Talagrand's recent proof of the Parisi formula, ultrametricity plays no apparent rôle, and it seems to be quite delicate to prove ultrametricity by Talagrand's method. This is quite curious as, on the other hand, ultrametricity plays a *crucial*

Received October 2004; revised July 2005.

[1]Supported in part by a grant from the Schweizerische Nationalfonds, contract no 20-100536/1.

rôle in the physicists nonrigorous derivation of the free energy, be that using the replica trick or the cavity method.]

We present here a simple and, as we think, natural generalization of the GREM which has *no* built in ultrametric structure. We, however, show that, in the limit, the model is ultrametrically organized. In this paper we address only the free energy. The more delicate investigation of the ultrametricity of the Gibbs distribution will be investigated in a forthcoming paper.

Throughout this paper, we fix a number $n \in \mathbb{N}$, and consider the set $I = \{1, \ldots, n\}$, as well as a collection of positive real numbers $\{a_J\}_{J \subset I}$ such that

$$\sum_{J \subset I} a_J = 1.$$

For convenience, we put $a_\varnothing \stackrel{\text{def}}{=} 0$. The relevant subset of $I$ will be only the ones with positive $a$-value. For $A \subset I$, we set

$$\mathcal{P}_A \stackrel{\text{def}}{=} \{J \subset A : a_J > 0\}, \qquad \mathcal{P} \stackrel{\text{def}}{=} \mathcal{P}_I.$$

For $n \in \mathbb{N}$, we set $\Sigma_N \stackrel{\text{def}}{=} \{1, \ldots, 2^N\}$. We also fix positive real numbers $\gamma_i$, $i \in I$, satisfying

$$\sum_{i=1}^n \gamma_i = 1,$$

and write $\Sigma_N^i \stackrel{\text{def}}{=} \Sigma_{\gamma_i N}$ where, for notational convenience, we assume that $2^{\gamma_i N}$ is an integer. For $N \in \mathbb{N}$, we will label the "spin configurations" $\sigma$ as

$$\sigma = (\sigma_1, \ldots, \sigma_n), \qquad \sigma_i \in \Sigma_N^i,$$

that is, we identify $\Sigma_N$ with $\Sigma_N^1 \times \cdots \times \Sigma_N^n$. For $J \subset I$, $J = \{j_1, \ldots, j_k\}$, $j_1 < j_2 < \cdots < j_k$, we write $\Sigma_{N,J} \stackrel{\text{def}}{=} \prod_{s=1}^k \Sigma_N^{j_s}$, and for $\sigma \in \Sigma_N$, we write $\sigma_J$ for the projected configuration $(\sigma_j)_{j \in J} \in \Sigma_{N,J}$. Our spin glass Hamiltonian is defined as

$$X_\sigma = \sum_{J \in \mathcal{P}} X^J_{\sigma_J}, \tag{1}$$

where $X^J_{\sigma_J}$, $J \in \mathcal{P}$, $\sigma_J \in \Sigma_{N,J}$ are independent centered Gaussian random variables with variance $a_J N$. The $X_\sigma$ are then Gaussian random variables with variance $N$ (Gaussian always means "centered Gaussian" through this note), but they are correlated. $\mathbb{E}$ will denote expectation with respect to these random variables. A special case is when $\mathcal{P} = \{I\}$, that is, when only $a_I \neq 0$, in which case it has to be one. Then the $X_\sigma$ are independent, that is, one considers simply a set of $2^N$ independent Gaussian random variables with variance $N$. This is the standard random energy model.

The generalized random energy model is a special case, too: It corresponds to the situation where the sets in $\mathcal{P}$ are "nested," meaning that $\mathcal{P}$ consists of an increasing sequence of subsets. Without loss of generality, we may assume that in this case

$$(2) \qquad \mathcal{P} = \{J_m : 1 \leq m \leq k\}, \qquad J_m \stackrel{\text{def}}{=} \{1, \ldots, n_m\},$$

where $1 \leq n_1 < n_2 < \cdots < n_k \leq n$. In the GREM case, the natural metric on $\Sigma_N$ coming from the covariance structure

$$d(\sigma, \sigma') \stackrel{\text{def}}{=} \sqrt{\mathbb{E}((X_\sigma - X_{\sigma'})^2)}$$

is an *ultrametric*. In the more general case (1) considered here, this metric is *not* an ultrametric.

To see this, take $n = 3$, $\mathcal{P} = \{\{1,2\}, \{1,3\}, \{2,3\}\}$, that is, where

$$(3) \qquad X_\sigma = X^{\{1,2\}}_{\sigma_1,\sigma_2} + X^{\{1,3\}}_{\sigma_1,\sigma_3} + X^{\{2,3\}}_{\sigma_2,\sigma_3},$$

with $a_J = 1/3$ for $J \in \mathcal{P}$. Then for $a, b, b', c, c' \in \Sigma_{N/3}$, $b \neq b'$, $c \neq c'$, one has

$$d((a,b,c),(a,b,c')) = d((a,b,c'),(a,b',c')) = \sqrt{2N/3},$$

whereas

$$d((a,b,c),(a,b',c')) = \sqrt{N},$$

contradicting ultrametricity.

Any of our models can be "coarse-grained" in many ways into a GREM. For that consider strictly increasing sequences of subsets of $I$: $\varnothing = A_0 \subset A_1 \subset \cdots \subset A_K = I$. We do not assume that the $A_i$ are in $\mathcal{P}$. We call such a sequence a *chain* $\mathbf{T} = (A_0, A_1, \ldots, A_K)$. We attach weights $\hat{a}_{A_j}$ to these sets by putting

$$(4) \qquad \hat{a}_{A_j} \stackrel{\text{def}}{=} \sum_{B \in \mathcal{P}_{A_j} \setminus \mathcal{P}_{A_{j-1}}} a_B.$$

Evidently $\sum_{j=1}^{K} \hat{a}_{A_j} = 1$, and if we assign random variables $X_\sigma(\mathbf{T})$, according to (1), we arrive after an irrelevant renumbering of $I$ at a GREM of the form (2). In particular, the corresponding metric $d$ is an ultrametric.

We write $\operatorname{tr}(\cdot)$ for averaging over $\Sigma_N$ (i.e., the coin-tossing expectation if we identify $\Sigma_N$ with $\{H,T\}^N$).

For a function $x : \Sigma_N \to \mathbb{R}$, set

$$Z_N(\beta, x) \stackrel{\text{def}}{=} \operatorname{tr} \exp[\beta x], \qquad F_N(\beta, x) \stackrel{\text{def}}{=} \frac{1}{N} \log(Z_N(\beta, x)),$$

and define the usual finite $N$ partition function, and free energy by

$$Z_N(\beta) \stackrel{\text{def}}{=} Z_N(\beta, X), \qquad F_N(\beta) \stackrel{\text{def}}{=} F_N(\beta, X), \qquad f_N(\beta) \stackrel{\text{def}}{=} \mathbb{E}(F_N(\beta, X)),$$

where $X$ is interpreted as random function $\Sigma_N \to \mathbb{R}$.

For any chain **T**, we attach to our model a GREM $(X_\sigma(\mathbf{T}))_{\sigma \in \Sigma_N}$, as explained above, and then

$$f_N(\mathbf{T}, \beta) \stackrel{\text{def}}{=} \mathbb{E}(F_N(\beta, X(\mathbf{T}))),$$

$$f(\mathbf{T}, \beta) \stackrel{\text{def}}{=} \lim_{N \to \infty} f_N(\mathbf{T}, \beta).$$

For a GREM, the limiting free energy is known to exist, and can be expressed explicitly, but in a somewhat complicated way (see [1, 3]). Our main result is that our generalization of the GREM does not lead to anything new in $N \to \infty$ limit, shedding hopefully some modest light on the "universality" of ultrametricity.

THEOREM 1.

(5) $$f(\beta) \stackrel{def}{=} \lim_{N \to \infty} f_N(\beta)$$

*exists and is also the almost sure limit of* $F_N(\beta)$.

$f(\beta)$ *is the free energy of a GREM. More precisely, there exists a chain* **T** *such that*

(6) $$f(\beta) = f(\mathbf{T}, \beta), \qquad \beta \geq 0.$$

$f(\mathbf{T}, \beta)$ *is minimal in the sense that*

(7) $$f(\beta) = \min_{\mathbf{S}} f(\mathbf{S}, \beta),$$

*the minimum being taken over all chains* **S**.

The fact that free energy is *self-averaging*, meaning that $f(\beta)$ (if the limit exists) is also the almost sure limit of the $F_N$, is a simple consequence of the Gaussian concentration inequality. We write $F_N$ as a function of the standardized variables $X^J_{\sigma_J}/\sqrt{a_J N}$. As

$$\left| \log \sum_i e^{a_i} - \log \sum_i e^{a'_i} \right| \leq \max_i |a_i - a'_i|, \qquad a_i, a'_i \in \mathbb{R},$$

we get that $F_N(\beta)$, regarded as a function of the collection $(X^J_{\sigma_J}/\sqrt{a_J N})$, is Lipschitz continuous with Lipschitz constant $\beta/\sqrt{N}$. By the usual concentration of measure estimates for Gaussian distributions (see, e.g., Proposition 2.18 of [6]), we have

(8) $$\mathbb{P}[|F_N(\beta) - \mathbb{E}F_N(\beta)| > \epsilon] \leq 2\exp\left[-\frac{\epsilon^2}{2\beta^2}N\right].$$

Using the Borel–Cantelli lemma, one sees that if $\lim_{N\to\infty} f_N(\beta)$ exists, then the $F_N(\beta)$ converge almost surely to this limit, too, and if $\lim_{N\to\infty} F_N(\beta)$ exists almost surely, then the limit is nonrandom and equals $\lim_{N\to\infty} f_N(\beta)$.

As for the strategy of the proof, the existence of the limit is established through a quite standard application of the *second moment method*, akin to that originally exploited by Derrida in his seminal paper [1]; this allows to express the limiting free energy in terms of a variational problem, which we then solve inductively.

For the reader's convenience, we briefly describe the mechanism which lies behind Theorem 1 for the Hamiltonian (3), but we allow for general (positive) variances $a_{12}, a_{13}, a_{23}$, and general $\gamma_i$. It is best to count the number of configurations $\sigma$ which reach a certain energy level $\lambda N$. It is evident that only an exponentially small portion of the total number $2^N$ of configurations achieve this, roughly formulated (we will be more precise later),

$$\#\{\sigma : X_\sigma \simeq \lambda N\} \simeq 2^N e^{-\rho(\lambda)N}, \qquad \rho(\lambda) > 0.$$

The free energy is obtained by the Legendre transform of $\rho$. In order to determine $\rho(\lambda)$, we count individually for each of the three parts in (3) how many configurations reach respective levels

(9)
$$\hat{\rho}(\lambda_1, \lambda_2, \lambda_3)$$
$$\simeq -\frac{1}{N} \log \#\{\sigma : X^{\{1,2\}}_{\sigma_1,\sigma_2} \simeq \lambda_1 N, \ X^{\{1,3\}}_{\sigma_1,\sigma_3} \simeq \lambda_2 N, \ X^{\{2,3\}}_{\sigma_2,\sigma_3} \simeq \lambda_3 N\} + \log 2,$$

with $\lambda_1 + \lambda_2 + \lambda_3 = \lambda$. Evidently,

(10) $$\rho(\lambda) = \inf_{\lambda_1+\lambda_2+\lambda_3=\lambda} \hat{\rho}(\lambda_1, \lambda_2, \lambda_3).$$

It turns out that one can get $\hat{\rho}$ by computing expectations inside the logarithm, provided only some naturally defined restrictions on the $\lambda_i$ are satisfied. For small $\lambda$, it is easily seen that one has an "equipartition" property, and that the optimal $\lambda_1, \lambda_2, \lambda_3$ are proportional to the respective variances, that is, $\lambda_1 = a_{12}\lambda$, $\lambda_2 = a_{13}\lambda$, $\lambda_3 = a_{23}\lambda$, and from that one obtains

(11) $$\rho(\lambda) = \lambda^2/2,$$

which is the same as if the $X_\sigma$ would be uncorrelated. Increasing $\lambda$, we, however, encounter restrictions from the structure of the Hamiltonian. First of all, $\lambda_1$ has to be such that there *are* any $\sigma_1, \sigma_2$ with $X^{\{1,2\}}_{\sigma_1,\sigma_2} \simeq \lambda_1 N = a_{12}\lambda N$. There are $2^{(\gamma_1+\gamma_2)N}$ pairs $(\sigma_1, \sigma_2)$, and as the $X^{\{1,2\}}_{\sigma_1,\sigma_2}$ are independent, the restriction is

$$2^{(\gamma_1+\gamma_2)N} \exp\left[-\frac{\lambda^2 a_{12} N}{2}\right] \gtrapprox 1.$$

(We are not considering any log-corrections.) This leads to the restriction

$$\lambda \leq \sqrt{\frac{2(\gamma_1 + \gamma_2)\log 2}{a_{12}}} \tag{12}$$

for the validity of (11), and there are two similar restrictions coming from $X^{\{1,3\}}$ and $X^{\{2,3\}}$. Even if these three restrictions are satisfied, it can be that there are simply totally not enough triples $(\sigma_1, \sigma_2, \sigma_3)$ left. A *necessary* condition for this is certainly that the expected number of $\#\{\sigma : X_\sigma \simeq \lambda N\}$ is not exponentially decaying, which is simply the condition that $\lambda \leq \sqrt{2\log 2}$. The somewhat astonishing fact is that these are the only conditions one has to take into considerations for the validity of (11). Now, there are two cases:

CASE 1. $\lambda \leq \sqrt{2\log 2}$ implies the other ones, that is,

$$\min_{1 \leq i < j \leq 3} \frac{\gamma_i + \gamma_j}{a_{ij}} \geq 1. \tag{13}$$

In that case, we are simply left with the restriction $\lambda \leq \sqrt{2\log 2}$, and the free energy is

$$f(\beta) = \sup_{\lambda \leq \sqrt{2\log 2}} (\beta\lambda - \lambda^2/2),$$

which is the free energy of an REM. In that case the internal structure of the model is irrelevant, at least for the free energy.

CASE 2. (13) is violated. For definiteness, assume that $(\gamma_1 + \gamma_2)/a_{12}$ is the smallest one.

In that case, (11) is only correct in the region (12). For $\lambda$ larger, there is no $(\sigma_1, \sigma_2)$ with $X^{\{1,2\}}_{\sigma_1,\sigma_2} \simeq a_{12}\lambda N$ (with probability close to 1), the maximum of the $X^{\{1,2\}}_{\sigma_1,\sigma_2}$ being at $m_{12}N$ ($\pm$ log-corrections), where

$$m_{12} \stackrel{\text{def}}{=} \sqrt{2(\gamma_1 + \gamma_2)a_{12}\log 2}.$$

Therefore, one has to restrict in (10) to $\lambda$'s with $\lambda_1 = m_{12}$. The only configurations $\sigma$ for which $X_\sigma \simeq \lambda N$ have to satisfy

$$X^{\{1,2\}}_{\sigma_1,\sigma_2} \simeq m_{12}N, \tag{14}$$

but there are now only subexponentially many $(\sigma_1, \sigma_2)$ left which achieve this feat, and the difference to $\lambda N$ has to be made by the field

$$Y_{\sigma_1,\sigma_2}(\sigma_3) \stackrel{\text{def}}{=} X^{\{1,3\}}_{\sigma_1,\sigma_3} + X^{\{2,3\}}_{\sigma_2,\sigma_3}, \qquad 1 \leq \sigma_3 \leq 2^{\gamma_3 N},$$

restricting $(\sigma_1, \sigma_2)$ to the few which satisfy (14). There is an upper limit $\lambda_{\max}$ for $\lambda$'s such that there are *any* $\sigma_3$ with $Y(\sigma_3) \simeq (\lambda - m_{12})N$. $\lambda_{\max} N - m_{12} N$ is simply the maximum of $2^{\gamma_3 N}$ independent Gaussians with variance $(a_{13} + a_{23})N$, that is,

$$\lambda_{\max} - m_{12} \stackrel{\text{def}}{=} \sqrt{2\gamma_3(a_{13} + a_{23})\log 2}.$$

The situation is similar to the one in the GREM with the only difference that, for $(\sigma_1, \sigma_2) \neq (\sigma_1', \sigma_2')$, the fields $Y_{\sigma_1,\sigma_2}$ and $Y_{\sigma_1',\sigma_2'}$ are not independent, except when $\sigma_1 \neq \sigma_1'$ and $\sigma_2 \neq \sigma_2'$. It is, however, fairly evident that, among the $(\sigma_1, \sigma_2)$ for which $X_{\sigma_1,\sigma_2}^{\{1,2\}} \simeq m_{12}N$ there will be no pairs with such a partial overlap, with probability close to 1, and therefore, it is quite natural one can handle the field $Y_{\sigma_1,\sigma_2}$ as if it would come from a second level of a two-level GREM. In fact, it turns out that, in the Case 2, the tree of Theorem 1 is $\{\{1,2\}, \{1,2,3\}\}$, and we replace our model with the coarse grained one with Hamiltonian $X'_{\alpha_1} + X''_{\alpha_1,\alpha_2}$, where $\#\alpha_1 = 2^{(\gamma_1 + \gamma_2)N}$, $\text{var}(X') = a_{12}N$, $\#\alpha_2 = 2^{\gamma_3 N}$, $\text{var}(X'') = (a_{13} + a_{23})N$.

This way of reasoning works for the general case. There are two issues which might be somewhat surprising. The first is that expressions (9) can always be evaluated by computing expectations inside the logarithm, provided one keeps some fairly trivial restrictions on the $\lambda_i$. Second, it is not entirely evident why these restrictions finally always lead to tree structures.

It is also interesting that the system always chooses from the many GREMs which can be obtained by coarse-grainings the one with minimal free energy. A similar behavior has already been obtained for the GREM itself in [2].

**2. Second moment estimates.** We fix some notation: If $(a_N)_{N \in \mathbb{N}}$ and $(b_N)_{N \in \mathbb{N}}$ are two sequences of positive real numbers, we write $a_N \asymp b_N$ if, for all $\varepsilon > 0$, there exists $N_0(\varepsilon) \in \mathbb{N}$ such that

$$e^{-\varepsilon N} b_N \leq a_N \leq e^{-\varepsilon N} b_N$$

for $N \geq N_0$. We also write $a_N \ll b_N$ if, for some $\delta > 0$, one has $a_N \leq b_N e^{-\delta N}$, again for large enough $N$. In that case, we also write $a_N = \Omega(b_N)$. The same notation are used in the case of sequences of random variables, just meaning that the relations hold almost surely (and therefore $N_0$ may depend on $\omega$). For $A \subset I$ (not necessarily in $\mathcal{P}$), we set

$$\gamma(A) \stackrel{\text{def}}{=} \sum_{i \in A} \gamma_i, \qquad \alpha(A) \stackrel{\text{def}}{=} \sum_{J \in \mathcal{P}_A} a_A.$$

We rewrite $F_N$ in terms of energy levels. For a collection $\lambda = (\lambda_J)_{J \in \mathcal{P}}, \lambda_J \in \mathbb{R}$ and $A \subset I$, we set

$$\mathcal{N}_{N,A}(\lambda) \stackrel{\text{def}}{=} \#\{\sigma \in \Sigma_{N,A} : X_{\sigma_J}^J \geq \lambda_J N, \forall J \in \mathcal{P}_A\},$$

$$\mathcal{N}_N(\lambda) \stackrel{\text{def}}{=} \mathcal{N}_{N,I}(\lambda).$$

Clearly,

(15) $$\{\mathcal{N}_{N,A}(\lambda) = 0\} \subset \{\mathcal{N}_N(\lambda) = 0\}.$$

We express $F_N$ in terms of the $\mathcal{N}_N(\lambda)$:

(16)
$$F_N(\beta) = \frac{1}{N} \log 2^{-N} (\beta N)^{|\mathcal{P}|} \int_{\mathbb{R}^{\mathcal{P}}} d\lambda \mathcal{N}_N(\lambda) \prod_{J \in \mathcal{P}} e^{\beta \lambda_J N}$$
$$= \frac{1}{N} \log \int_{\mathbb{R}^{\mathcal{P}}} d\lambda \mathcal{N}_N(\lambda) \prod_{J \in \mathcal{P}} e^{\beta \lambda_J N} - \log 2 + O\left(\frac{\log N}{N}\right).$$

We first want to take out the $\lambda$ for which $\mathcal{N}_N(\lambda) = 0$ for large $N$. As these are integer valued random variables, it is clear that $\mathbb{E}\mathcal{N}_N(\lambda) \ll 1$ implies $\mathcal{N}_N(\lambda) = 0$ for large enough $N$, almost surely. It, however, turns out that this condition is not sufficient for our purpose, but remark that, if for *some* $A \subset I$, one has $\mathbb{E}\mathcal{N}_{N,A}(\lambda) \ll 1$, then, by (15), one has $\mathcal{N}_N(\lambda) = 0$ for large enough $N$ as well.

LEMMA 2. (a) *For any $\lambda \in \mathbb{R}^{\mathcal{P}}$ and $A \subset I$, we have*

$$\mathbb{E}\mathcal{N}_{N,A}(\lambda) \asymp 2^{\gamma(A)N} \exp\left[-\sum_{J \in \mathcal{P}_A} \frac{(\lambda_J^+)^2}{2 a_J} N\right],$$

*where $\lambda_J^+ \stackrel{\text{def}}{=} \max(\lambda_J, 0)$.*
  (b) *There exists $C > 0$ such that*

$$\mathbb{E}\mathcal{N}_N(\lambda) \leq C 2^N \exp\left[-\sum_{J \in \mathcal{P}} \frac{(\lambda_J^+)^2}{2 a_J} N\right]$$

*for all $\lambda \in \mathbb{R}^{\mathcal{P}}$, and all $N$.*
  (c) *Let $\lambda \in \mathbb{R}^{\mathcal{P}}$. If for some $A \subset I$ one has*

$$\sum_{J \in \mathcal{P}_A} \frac{(\lambda_J^+)^2}{2 a_J} > \gamma(A) \log 2,$$

*then $\mathbb{P}(\mathcal{N}_N(\lambda) \neq 0) \ll 1$ and, in particular, $\mathcal{N}_N(\lambda) = 0$ for large enough $N$, almost surely.*

PROOF. (a) and (b) follow by standard Gaussian tail estimates, and in case (c), by (15), we have

$$\mathbb{P}(\mathcal{N}_N(\lambda) \neq 0) \leq \mathbb{P}(\mathcal{N}_{N,A}(\lambda) \neq 0) \leq \mathbb{E}\mathcal{N}_{N,A}(\lambda),$$

which proves (c). □

Let

$$\Delta \stackrel{\text{def}}{=} \left\{ \lambda \in \mathbb{R}^{\mathcal{P}} : \sum_{J \in \mathcal{P}_A} \frac{(\lambda_J^+)^2}{2a_J} \leq \gamma(A) \log 2, \forall A \subset I \right\},$$

$$\Delta^+ \stackrel{\text{def}}{=} \{\lambda \in \Delta : \lambda_J \geq 0, \forall J \in \mathcal{P}\}.$$

LEMMA 3. *If $\lambda \in \operatorname{int} \Delta$, then*

(17) $$\mathcal{N}_N(\lambda) \asymp \mathbb{E}\mathcal{N}_N(\lambda) \asymp \exp\left[N\left(\log 2 - \sum_{J \in \mathcal{P}} (\lambda_J^+)^2 / 2a_J\right)\right].$$

PROOF. The second relation is Lemma 2(a). For the proof of the first, it suffices to show that $\lambda \in \operatorname{int} \Delta$ implies

(18) $$\operatorname{var} \mathcal{N}_N(\lambda) \ll (\mathbb{E}\mathcal{N}_N(\lambda))^2.$$

In fact, from (18), Chebyshev's inequality and the Borel–Cantelli lemma immediately imply (17).

We abbreviate $\mathbb{P}(X_{\sigma_J}^J \geq N\lambda_J)$ by $p_J(N)$ ($\lambda$ is kept fixed through this proof). With this notation, $\mathbb{E}\mathcal{N}_N(\lambda) = 2^N \prod_{J \in \mathcal{P}} p_J(N)$:

$$\begin{aligned}
\mathbb{E}\mathcal{N}_N(\lambda)^2 &= \sum_{\sigma,\sigma' \in \Sigma_N} \prod_{J \in \mathcal{P}} \mathbb{P}(X_{\sigma_J}^J \geq N\lambda_J, \ X_{\sigma'_J}^J \geq N\lambda_J) \\
&= \sum_{A \subset I} \sum_{(\sigma,\sigma') \in \Lambda_A} \prod_{J \in \mathcal{P}} \mathbb{P}(X_{\sigma_J}^J \geq N\lambda_J, \ X_{\sigma'_J}^J \geq N\lambda_J) \\
&= \sum_{A \subset I} |\Lambda_A(N)| \prod_{J \in \mathcal{P}_A} p_J(N) \prod_{J \in \mathcal{P} \setminus \mathcal{P}_A} p_J(N)^2,
\end{aligned}$$

where $\Lambda_A(N)$ consists of those pairs $(\sigma, \sigma')$ which agree on $A$ and disagree on $I \setminus A$. For $A = \varnothing$, $2^{2N} - |\Lambda_\varnothing(N)| \ll 2^{2N}$, and therefore,

$$\operatorname{var} \mathcal{N}_N(\lambda) = \sum_{A \neq \varnothing} |\Lambda_A(N)| \prod_{J \in \mathcal{P}_A} p_J(N) \prod_{J \in \mathcal{P} \setminus \mathcal{P}_A} p_J(N)^2 + \Omega((\mathbb{E}\mathcal{N}_N(\lambda))^2).$$

$$|\Lambda_A(N)| = 2^{\gamma(A)N} \prod_{i \notin A} 2^{\gamma_i N}(2^{\gamma_i N} - 1) = 2^{2N} 2^{-\gamma(A)N} + \Omega(|\Lambda_A(N)|).$$

As by assumption,

$$2^{-\gamma(A)N} \ll \prod_{J \in \mathcal{P}_A} p_J(N) \asymp \exp\left[-\sum_{J \in \mathcal{P}_A} \frac{(\lambda_J^+)^2}{2a_J} N\right], \qquad A \neq \varnothing,$$

we have, for any $A \neq \varnothing$,

$$|\Lambda_A(N)| \prod_{J \in \mathcal{P}_A} p_J(N) \prod_{J \in \mathcal{P} \setminus \mathcal{P}_A} p_J(N)^2 \ll 2^{2N} \prod_{J \in \mathcal{P}} p_J(N)^2,$$

proving $\operatorname{var} \mathcal{N}_N(\lambda) \ll (\mathbb{E} \mathcal{N}_N(\lambda))^2$. □

Let

$$\psi(\lambda, \beta) \stackrel{\text{def}}{=} \sum_{J \in \mathcal{P}} \left( \beta \lambda_J - \frac{\lambda_J^2}{2 a_J} \right). \tag{19}$$

PROPOSITION 4. *The free energy as defined in* (5) *exists and is given as*

$$f(\beta) = \sup_{\lambda \in \Delta^+} \psi(\lambda, \beta). \tag{20}$$

PROOF. We show the lower bound for $\liminf_{N \to \infty} F_N(\beta)$ and the upper bound for $\limsup_{N \to \infty} f_N(\beta)$. By the self-averaging property (8), this proves the statement.

We use the integral representation (16). If $\mu = (\mu_J)$, $\nu = (\nu_J)$ satisfy $\mu_J < \nu_J$ for all $J$, we write

$$[\mu, \nu) \stackrel{\text{def}}{=} \{\lambda : \mu_J \leq \lambda_J < \nu_J, \forall J \in \mathcal{P}\}.$$

If $[\mu, \nu) \subset \Delta^+$, we have

$$\liminf_{N \to \infty} F_N(\beta) \geq \liminf_{N \to \infty} \frac{1}{N} \log \int_{[\mu,\nu)} d\lambda \, \mathcal{N}_N(\lambda) \prod_{J \in \mathcal{P}} e^{\beta \lambda_J N}$$

$$\geq \liminf_{N \to \infty} \frac{1}{N} \log \mathcal{N}_N(\nu) \int_{[\mu,\nu)} d\lambda \prod_{J \in \mathcal{P}} e^{\beta \lambda_J N}$$

$$\geq \sum_J \left( \beta \mu_J - \frac{\nu_J^2}{2 a_J} \right).$$

As this holds for arbitrary $[\mu, \nu) \subset \Delta^+$, $\liminf_{N \to \infty} F_N(\beta) \geq \sup_{\lambda \in \Delta^+} \psi(\lambda, \beta)$ follows.

For the upper bound, let $\Delta_\varepsilon$ be an $\varepsilon$-neighborhood of $\Delta$. From Lemma 2(c), we have

$$\mathbb{P}\left( F_N(\beta) \neq \frac{1}{N} \log \int_{\Delta_\varepsilon} d\lambda \, \mathcal{N}_N(\lambda) e^{\beta \sum_J \lambda_J N} \right) \ll 1,$$

and therefore,

$$\left| \mathbb{E} F_N(\beta) - \frac{1}{N} \mathbb{E} \log \int_{\Delta_\varepsilon} d\lambda \, \mathcal{N}_N(\lambda) e^{\beta \sum_J \lambda_J N} \right| \ll 1.$$

By Jensen's inequality and Lemma 2(b), we have

$$\limsup_{N \to \infty} \frac{1}{N} \mathbb{E} \log \int_{\Delta_\varepsilon} d\lambda \, \mathcal{N}_N(\lambda) e^{\beta \sum_J \lambda_J N}$$

$$\leq \limsup_{N \to \infty} \frac{1}{N} \log \int_{\Delta_\varepsilon} d\lambda \, \mathbb{E} \mathcal{N}_N(\lambda) e^{\beta \sum_J \lambda_J N}$$

$$\leq \limsup_{N \to \infty} \frac{1}{N} \log 2^N \int_{\Delta_\varepsilon} d\lambda \exp\left[ N \sum_{J \in \mathcal{P}} \left( \beta \lambda_J - \frac{(\lambda_J^+)^2}{2a_J} \right) \right]$$

$$\leq \sup_{\lambda \in \Delta_\varepsilon} \psi(\lambda, \beta).$$

As $\varepsilon > 0$ is arbitrary, $\limsup_{N \to \infty} \mathbb{E} F_N(\beta) \leq \sup_{\lambda \in \Delta^+} \psi(\lambda, \beta)$ follows. $\square$

**3. The optimization problem.** We first discuss the special case of a GREM. Therefore, we assume that the sets in $\mathcal{P}$ are nested, that is, $\mathcal{P} = \{J_1, \ldots, J_m\}$, where $\varnothing \subset J_1 \subset \cdots \subset J_m$. If $A \subset I$, put $l_A \stackrel{\text{def}}{=} \max\{l : J_l \subset A\}$. Evidently,

$$\sum_{J \in \mathcal{P}_A} \frac{\lambda_J^2}{2a_J} \leq \log 2 \, \gamma(A)$$

follows from

$$\sum_{i=1}^{l_A} \frac{\lambda_{J_i}^2}{2a_{J_i}} \leq \log 2 \, \gamma(J_{l_A}).$$

Therefore, $\lambda \in \Delta^+$ is equivalent with

$$\sum_{i=1}^{l} \frac{\lambda_{J_i}^2}{2a_{J_i}} \leq \log 2 \, \gamma(J_l), \qquad 1 \leq l \leq m,$$

(and, of course, that all components are nonnegative). Therefore, we have proved the following:

LEMMA 5. *Assume that $\mathcal{P}$ is nested as above. Then*

$$f(\beta) = \sup \left\{ \psi(\lambda, \beta) : \sum_{i=1}^{l} \frac{\lambda_{J_i}^2}{2a_{J_i}} \leq \log 2 \, \gamma(J_l), \, 1 \leq l \leq m \right\}.$$

This lemma proves that, in our more general situation, for any chain $\mathbf{T}$, the corresponding GREM free energy is an upper bound.

COROLLARY 6. *For any chain $\mathbf{T}$, we have*

$$f(\beta) \leq f(\mathbf{T}, \beta), \qquad \beta \geq 0.$$

PROOF. For a given chain $\varnothing = A_0 \subset A_1 \subset A_2 \subset \cdots \subset A_K = I$, we consider $\Delta_{\mathbf{T}}^+$ which is obtained by dropping the conditions for the $A$'s which are not in the chain. Then

$$f(\beta) = \sup_{\lambda \in \Delta^+} \psi(\lambda, \beta) \leq \sup_{\lambda \in \Delta_{\mathbf{T}}^+} \psi(\lambda, \beta).$$

We claim that

$$\sup_{\lambda \in \Delta_{\mathbf{T}}^+} \psi(\lambda, \beta) = f(\mathbf{T}, \beta),$$

which proves the corollary. To see this equation, we write

$$\psi(\lambda, \beta) = \sum_{j=1}^{K} \sum_{J \in \mathcal{P}_{A_j} \setminus \mathcal{P}_{A_{j-1}}} \left( \beta \lambda_J - \frac{\lambda_J^2}{2 a_J} \right) \stackrel{\text{def}}{=} \sum_{j=1}^{K} \psi_j(\lambda_j, \beta), \qquad \text{say,}$$

where $\lambda_j \stackrel{\text{def}}{=} (\lambda_J)_{J \in \mathcal{P}_{A_j} \setminus \mathcal{P}_{A_{j-1}}}$. Set

$$f_j(\beta, t) \stackrel{\text{def}}{=} \sup \left\{ \psi_j(\lambda_j, \beta) : \sum_{J \in \mathcal{P}_{A_j} \setminus \mathcal{P}_{A_{j-1}}} \frac{\lambda_J^2}{2 a_J} = t \right\} = \beta \sqrt{2 t \hat{a}_j} - t,$$

where $\hat{a}_j \stackrel{\text{def}}{=} \sum_{J \in \mathcal{P}_{A_j} \setminus \mathcal{P}_{A_{j-1}}} a_J$, that is, $f_j(\beta, s^2/2\hat{a}_j) = \beta s - s^2/2\hat{a}_j$. We therefore see that

$$\sup_{\lambda \in \Delta_{\mathbf{T}}^+} \psi(\lambda, \beta) = \sup \left\{ \sum_{j=1}^{K} \left( \beta s_j - \frac{s_j^2}{2\hat{a}_j} \right) : \sum_{j=1}^{l} \frac{s_j^2}{2\hat{a}_j} \leq \log 2 \, \gamma(A_l), \ 1 \leq l \leq K \right\}$$

$$= f(\mathbf{T}, \beta),$$

the last equality by Lemma 5. □

In order to finish the proof of Theorem 1, it only remains to construct a chain $\mathbf{T}$ which satisfies $f(\beta) \geq f(\mathbf{T}, \beta)$. Then one has also equality by the above corollary.

For $B \subset I$, let

$$\alpha(B) \stackrel{\text{def}}{=} \sum_{J \in \mathcal{P}_B} a_J$$

and for $B \subset A$, set

$$\rho(B, A) \stackrel{\text{def}}{=} \sqrt{\frac{2 \log 2 (\gamma(A) - \gamma(B))}{\alpha(A) - \alpha(B)}},$$

$$\hat{\rho}(B) \stackrel{\text{def}}{=} \min_{A : A \supset B, A \neq B} \rho(B, A).$$

We construct a strictly increasing sequence of subsets $A_0 \stackrel{\text{def}}{=} \varnothing \subset A_1 \subset A_2 \subset \cdots \subset A_K = I$ and parameters $\beta_0 \stackrel{\text{def}}{=} 0 < \beta_1 < \beta_2 < \cdots < \beta_K < \infty$ by recursion. Assume that $\varnothing \subset A_1 \subset A_2 \subset \cdots \subset A_k$ and $0 < \beta_1 < \beta_2 < \cdots < \beta_k$ are constructed such that the following conditions are satisfied:

C1($k$) $\beta_j = \hat{\rho}(A_{j-1})$, $j \leq k$.
C2($k$) For $j \leq k$ and any $A \supset A_{j-1}$ which satisfies $\beta_j = \rho(A_{j-1}, A)$, one has $A \subset A_j$, that is, $A_j$ is *maximal* with $\beta_j = \rho(A_{j-1}, A_j)$.

For $k = 0$, the conditions are void. If $A_k = I$, then the construction is finished, and we have $K \stackrel{\text{def}}{=} k$. Therefore, assume $A_k \neq I$. Then we set $\beta_{k+1} \stackrel{\text{def}}{=} \hat{\rho}(A_k)$, and prove first that $\beta_{k+1} > \beta_k$. We claim that, for any $A \supset A_k$, $A \neq A_k$, one has

$$2 \log 2(\gamma(A) - \gamma(A_k)) > \beta_k^2(\alpha(A) - \alpha(A_k)).$$

Indeed, because of $2 \log 2(\gamma(A_k) - \gamma(A_{k-1})) = \beta_k^2(\alpha(A_k) - \alpha(A_{k-1}))$,

$$2 \log 2(\gamma(A) - \gamma(A_k)) < \beta_k^2(\alpha(A) - \alpha(A_k))$$

would contradict condition C1($k$) and equality would contradict C2($k$).

It only remains to construct $A_{k+1}$ which satisfies C2($k+1$). Assume there are two sets $A, A' \supset A_k$, $A, A' \neq A_k$ satisfying

(21) $$\rho(A_k, A) = \rho(A_k, A') = \beta_{k+1}.$$

We claim that then also $\rho(A_k, A \cup A') = \beta_{k+1}$. Remark that

$$\alpha(A \cup A') \geq \alpha(A) + \alpha(A') - \alpha(A \cap A'),$$
$$\gamma(A \cup A') = \gamma(A) + \gamma(A') - \gamma(A \cap A'),$$

and therefore,

$$2 \log 2(\gamma(A \cup A') - \gamma(A_k)) - \beta_{k+1}^2(\alpha(A \cup A') - \alpha(A_k))$$
$$\leq 2 \log 2[\gamma(A) + \gamma(A') - \gamma(A \cap A') - \gamma(A_k)]$$
$$\quad - \beta_{k+1}^2[\alpha(A) + \alpha(A') - \alpha(A \cap A') - \alpha(A_k)]$$
$$= \beta_{k+1}^2[\alpha(A \cap A') - \alpha(A_k)] - 2 \log 2[\gamma(A \cap A') - \gamma(A_k)] \leq 0,$$

the equality by (21), and the last inequality by the definition of $\beta_{k+1}$. From the definition of $\beta_{k+1}$, we therefore conclude that

$$2 \log 2(\gamma(A \cup A') - \gamma(A_k)) = \beta_{k+1}^2(\alpha(A \cup A') - \alpha(A_k)).$$

We therefore find a *unique* maximal set $A_{k+1} \supset A_k$ which satisfies $\rho(A_k, A_{k+1}) = \beta_{k+1}$, and so we have constructed $\beta_{k+1} > \beta_k$, $A_{k+1} \supset A_k$, $A_{k+1} \neq A_k$ such that C1($k+1$) and C2($k+1$) are satisfied. The construction terminates after a finite number of steps.

We claim now that, with $\mathbf{T} \stackrel{\text{def}}{=} (\varnothing, A_1, \ldots, A_{K-1}, I)$, we have

(22) $$f(\beta) \geq f(\beta, \mathbf{T}).$$

Clearly, if $\beta > 0$ is small enough, the maximum in (20) is attained in $\lambda_J^{(1)}(\beta) \stackrel{\text{def}}{=} a_J \beta$ for all $J$, and therefore,

$$f(\beta) = \frac{\beta^2}{2}$$

for small $\beta$. This remains valid as long as $(\beta^2/2)\alpha(A) \leq \gamma(A)$ for all $A$, that is, for $\beta \leq \beta_1$. For $\beta_k < \beta \leq \beta_{k+1}$, we choose $\lambda^{(k+1)}(\beta)$ defined by

$$(23) \qquad \lambda_J^{(k+1)}(\beta) \stackrel{\text{def}}{=} \begin{cases} a_J \beta_m, & \text{for } J \in \mathcal{P}_{A_m} \setminus \mathcal{P}_{A_{m-1}},\ 1 \leq m \leq k, \\ a_J \beta, & \text{for } J \notin \mathcal{P}_{A_k}. \end{cases}$$

This choice (23) satisfies the side conditions in the range of $\beta$ we are considering, and hence,

$$(24) \qquad \psi(\lambda^{(k+1)}(\beta), \beta) \leq f(\beta).$$

We show now that $f(\beta, \mathbf{T}) = \psi(\lambda^{(k+1)}(\beta), \beta)$ for $\beta_k \leq \beta \leq \beta_{k+1}$. An elementary computation gives

$$\psi(\lambda^{(k+1)}(\beta), \beta) = \beta \sum_{i=1}^{k} \beta_i [\alpha(A_i) - \alpha(A_{i-1})] - \gamma(A_k) \log 2 + \frac{\beta^2}{2}(1 - \alpha(A_k))$$

$$= \beta \sum_{i=1}^{k} \beta_i \hat{a}(A_i) - \gamma(A_k) \log 2 + \frac{\beta^2}{2} \sum_{i=k+1}^{K} \hat{a}(A_i),$$

where $\hat{a}(A_i)$ is defined by (4). This is exactly the free energy of the corresponding GREM as given in [1]. [It is, in fact, elementary to check that $\lambda^{(k+1)}(\beta)$ is the maximizing vector $\lambda$ for the GREM corresponding to the above chain when $\beta_k \leq \beta \leq \beta_{k+1}$.] We have therefore proved Theorem 1.

REMARK 7. We have, in fact, proved that

$$f(\beta) = \beta \sum_{i=1}^{k} \beta_i [\alpha(A_i) - \alpha(A_{i-1})] - \gamma(A_k) \log 2 + \frac{\beta^2}{2}(1 - \alpha(A_k))$$

for $\beta_k \leq \beta \leq \beta_{k+1}$.

**Acknowledgment.** We wish to thank the referee for careful reading and valuable comments.

INSTITUT FÜR MATHEMATIK
UNIVERSITÄT ZÜRICH
WINTERTHURERSTRASEE 190
CH-8057 ZÜRICH
SWITZERLAND
E-MAIL: nkistler@amath.unizh.ch
        eb@math.unizh.ch